\documentclass{amsart}

\usepackage[pdftex, colorlinks=true, linkcolor=blue, urlcolor=blue]{hyperref}
\usepackage[usenames,dvipsnames,pdftex]{xcolor}

\usepackage{pgf}
\usepackage{xxcolor}
\usepackage{hyperref}

\usepackage{amsmath,amssymb,amsbsy, amsthm} 

\usepackage{xcolor}
\usepackage{fancyhdr, blindtext}

\issueinfo{8}{1}{number 1}{2014}

\copyrightinfo{2014}{\href{http://albanian-j-math.com}{albanian-j-math.com}}

\def\issn{{\sc \textbf{ISSN: }} 1930-1235; }
\def\issueyear{\textbf{2014}}

\PII{\issn  (\issueyear)}

\pagespan{3}{8}

\title{Equations for superelliptic curves over their minimal field of definition}

\keywords{superelliptic curves \and field of moduli \and minimal field of definition \and  Shaska invariants}
\subjclass[2000]{11G30 \and 11G50 \and 14G40}

\usepackage{amsrefs}

\newtheorem{thm}{Theorem}

\def\C{\mathbb C}

\def\C{{\mathbb C}}

\def\M{\mathcal M}
\def\X{\mathcal X}

\def\s{\mathfrak s}

\def\u{\mathfrak s}

\def\X{\mathcal X}

\def\d{{\delta }}

\def\D{\Delta}

\def\G{G}

\begin{document}
 
\maketitle

\hrule

\vspace{.6cm}



\begin{center}
\textsc{Lubjana Beshaj} \\

\smallskip

\begin{small}
\textit{
Department of Mathematics and Statistics\\
Oakland University\\
Rochester, MI, 48386. \\
Email:  beshaj@oakland.edu \\
}
\end{small}
\end{center}


\bigskip

\begin{center}
\textsc{Fred Thompson} \\

\smallskip

\begin{small}
\textit{
Department of Mathematics and Statistics\\
Oakland University\\
Rochester, MI, 48386. \\
Email:  fjthomps@oakland.edu \\
}
\end{small}
\end{center}


\vspace{.5cm}

\hrule

\vspace{.2cm}

\begin{abstract}
Let $\X_g$ be a genus $g\geq 2$ superelliptic curve, $F$ its field of moduli, and $K$ the minimal field of definition.  
In this short note we construct an equation of the curve $\X_g$ over its minimal field of definition $K$
when $\X_g$ has   extra automorphisms. We make use of the dihedral invariants of superelliptic curves as defined by Shaska in \cite{issac-03} and results on the automorphism groups of superelliptic curves as in \cite{Sa2}. 
\end{abstract}


\vspace{0.6cm}  \hrule    \vspace{.5cm}


\section{Introduction}


Given an algebraic curve $\X$ of genus $g\geq 2$, it is an open problem to determine an equation for $\X$ over its minimal field of definition $K$.  It is well known that the minimal field of definition is an algebraic extension of the field of moduli $F$. While for small genus it is known how to construct such equations, in general this is still an open problem.  The overall strategy is to describe the point in the moduli space $\M_g$ corresponding to the given curve.  This determines the field of moduli $F$ and the minimal field of definition $K$ is a finite extension of the field of moduli. However, describing the moduli point explicitly can be done only for superelliptic curves of small genus; see \cite{s1, s2, s3}.  

Superelliptic curves are curves with affine equation $y^n = f(x)$.  Such curves have at least an automorphism of order $n$.  The quotient by the automorphism group of such curves is a genus 0 curve, hence a conic.  This conic always has a rational point over a quadratic extension of the field of moduli.  Hence, for superelliptic curves $[K : F ] \leq 2$.  If the automorphism group of $\X$ is isomorphic to the cyclic group of order $n$ then an idea of Clebsch can be extended to determine if the field of moduli is a field of definition.  Moreover an equation can be determined over the minimal field of definition.  This is intended in \cite{bst}.

When the superelliptic curves have extra automorphisms, i.e. the automorphism group has size $ > n$ then the algorithm suggested above does not work.  The isomorphism classes of such curves are determined by dihedral invariants (or Shaska invariants)  as in \cite{g-sh, AK,  s4}. 

In this short note we give an equation of superelliptic curves  of genus $g\geq 2$ with extra automorphisms over the minimal field of definition $K$ and determine the algebraic conditions in terms of such invariants of  curves when the field of moduli is a field of definition. 

Our main result is the following.  Let $\X$ be a genus $g\geq 2$ superelliptic curve, defined over $\C$, with an extra automorphism, $\s_1, \dots , \s_g$  its dihedral invariants, $F$ the field of moduli, and $K$ its minimal field of definition.  Then,  \\
 
i) The minimal field of definition $K$  is    $K = F ( \sqrt{\Delta_\u})$   

ii)  The equation of $\X$ over $K$ is 
\[
y^n = A\, x^{\d (s+1)} + A \,  x^{\d s} + \sum_{i=1}^{s-1}  2^{s-i} \, \u_1 \cdot \frac  { \u_s^i \u_i - A \u_{s+1-i}   }    { 2^s \u_1^2 - \u_s^{s+1}  } \cdot x^{\d \cdot i} + 1  
\]
where \[ 2^{s+1} A^2 - 2^{s+1} \u_1 A + \u_s^{s+1} = 0. \]  
and $\D_\u$ is the discriminant of the above quadratic, 
\[ \D_\u =  2^{s+1} \left(  2^{s+1}  \u_1^2 - 4 \u_s^{s+1}  \right).   \]

Hence, this provides an Weierstrass equation of the curve over $k ( \sqrt{\D_\u} )$. 

An immediate consequence of the above result is that the field of moduli is a field of definition when the above quadratic has rational solutions.  This happens if and only if $\D_s$ is a complete square.  

It was noted in \cite{g-sh} that when $\D_s$ the automorphism group of the curve is larger and can be explicitly determined. The case when the genus $g$ is odd differs from the case when it is even. As a corollary we get that if $\D_s =0$ then the field of moduli is a field of definition as noted in  \cite{g-sh} for hyperelliptic curves. 

The results of this paper determine when the field of moduli is a field of definition and give an equation of the curve over the minimal field of definition for almost all superelliptic curves with extra automorphism.  The next natural thing to study is the case of the generic superelliptic curve, that is the curves with equation $y^n=f(x)$ and automorphism group of order $n>2$.  Such algorithm is given in \cite{Me} for genus $g=2$ and it is intended in \cite{bst} for all superelliptic curves. 


\section{Preliminaries} 
Let $\X_g$ be a genus $g \geq 2$ curve with full automorphism group $G=Aut(\X_g)$.  The curve $\X_g$ is called a \textbf{superelliptic curve} if there exists an element $\tau \in G$ which is central in $G$ and $g(\X_g/\langle \tau \rangle)=0$.  Denote by $H$ the cyclic group generated by $\tau$,  $H= \langle \tau \rangle$.  Thus, $\overline G =G/H$ is called the reduced automorphism group of $\X_g$ with respect to $H$.  

Superelliptic curves are curves with affine equation $y^n=f(x)$. Denote with $K=k(x, y)$ the function field of $\X_g$ and by $k(x)$ the genus zero subfield of $K$ fixed by $H$.  Then, $[K:k(x)]=n$, where $n=|H|$.  The group $\overline G$ is a subgroup  of the group of automorphisms of a genus zero curve. Therefore,  $\overline G < PGL_2(k)$ and $\overline G$ is finite. Then,  $\overline G$ is isomorphic to one of the following groups $C_m, D_m, A_4, S_4, A_5$. Since $G$ is a degree $n$ extension of $\overline G$ and we know the possible groups that occur as $\overline G$, it is possible to determine $G$ and the equation of $K$, see \cite{Sa2}. 

The group $\G$ acts on $k(x)$ via the natural way. The fixed field of this action is a genus 0 field, say $k(z)$. Thus, $z$
is a degree $|\G|$ rational function in $x$, say $z=\phi(x)$.

Given a superelliptic curve $\X_g$  with equation $y^n=f(x)$ such that $\Delta (f,x) \neq 0$,  the genus of the curve can be calculated using the following formula
\[ g = 1 + \frac 1 2 \, \left( nd -n -d - \gcd (d, n) \right). \]
where $\deg f =d >n$. If $d$ and $n$ are relatively prime then $g = \frac {(n-1)(d-1)} 2$, see \cite{tw2} for proof. 

Much interesting to us are superelliptic curves with extra automorphism. Let 
$\X_g$ be a superelliptic curve that has an extra automorphism $\sigma \in G$ such that its projection $\overline \sigma \in \overline G$ has order $\delta \geq 2$. Then the equation of the superelliptic curve is given as $y^n=g(x^\delta)$ or $y^n=xg(x^\delta)$, for some $g \in k[x]$, see \cite{s4} for proof. 
In other words $\X_g$ has equation
\[y^n= g(x^\delta):= x^{s \delta }+a_{s-1}x^{(s -1) \delta}+ \cdots+a_1x^\delta +1,\]
or
\[y^n= xg(x^\delta):= x^{(s+1) \delta }+a_sx^{s  \delta}+ \cdots+a_1x^\delta +x.\]
For both cases the dihedral invariants of such curves or \textit{Shaska-invariants} denoted by $\s$-invariants  are defined   in \cite{issac-03, AK, g-sh}. They were discovered by Shaska in his thesis for curves of genus 2 with extra automorphisms and later generalized to all hyperelliptic curves in \cite{g-sh} for all hyperelliptic curves with extra automorphisms.  Such invariants are used by many authors in computational aspects of hyperelliptic and superelliptic curves such as Duursma, Ritzenthaler, Lauter, Lercier, et al.
We define them for our purposes in the next section.


\section{Equation of superelliptic curves and dihedral invariants}

Let $\X_g$ be a superelliptic curve defined over a field $k$, $\mbox{char} k =0$ such that $\X_g$ has an extra involution and its Weierstrass equation   is given by 
\begin{equation}\label{w-eq} y^n= x^{\d (s+1)} + a_s x^{\d s} + a_{s-1} x^{\d (s-1)} + \cdots + a_2 x^{\d \cdot 2} + a_1 x^{\d} + 1    \end{equation}
Our main goal is to find an equation of this curve defined over its minimal field of definition.   The corresponding moduli point of such curves is determined by the dihedral invariants $\u_1, \dots , \u_s$ and the field of moduli is $k (\u_1, \dots , \u_s)$. 

Recall that the dihedral invariants are defined as follows
\[
\begin{split}
& \u_1  = a_1^{s+1} + a_s^{s+1}  \\
& \u_2  = a_1^{s-1} a_2 + a_s^{s-1}a_{s-1} \\
& \cdots \\
& \u_i  = a_1^{s+1-i} a_i + a_s^{s+1-i}a_{s+1-i} \\
& \cdots \\
& \u_{s+1-i} = a_1^{i} a_{s+1-i} + a_s^{i} a_{i}\\
 & \cdots \\
&  \u_{s-1} = a_1^2 a_{s-1} + a_s^2 a_2 \\
& \u_s = 2 a_1 a_s \\ 
\end{split}
\]
Notice that these invariants are homogenous polynomials of degree $s+1$ to 2 respectively.  The field of moduli of the corresponding curve is given by $k(\u_1, \dots , \u_s)$.  Our goal is to find a Weierstrass equation over $k(\u_1, \dots , \u_s)$ of the curve in Eq.~\eqref{w-eq}.

We perform a coordinate change \[x \to \sqrt[\d]{ a_s} \, x\]
to get 
\[ y^n = a_s^{s+1} \, x^{\d (s+1)} + a_s^{s+1} \,  x^{\d s} + a_{s-1} \cdot a_s^{s-1} \,  x^{\d (s-1)} + \cdots + a_2 \cdot a_s^2 \, x^{\d \cdot 2} + a_1 a_s \, x^{\d} + 1    \]
Denote by $A:=a_s^{s+1}$.   Then we have 
\[ 2^{s+1} A^2 - 2^{s+1} \u_1 A + \u_s^{s+1} = 0. \]  
This quadratic has discriminant 
\[\Delta_\u=  2^{s+1} \left(  2^{s+1}  \u_1^2 - 4 \u_s^{s+1}  \right) \]
The equation of the curve becomes 
\[ y^n =  A\, x^{\d (s+1)} + A \,  x^{\d s} + \sum_{i=1}^{s-1} a_i  a_s^i  \cdot x^{\d \cdot i} + 1  \]
We will show that all coefficients $a_i a_s^i$, $i=1, \dots, s-1,$ can be expressed in terms of the dihedral invariants and $A$.  Hence, we have an equation of the curve over the quadratic extension $k \sqrt{ \Delta_\u}$.

\begin{thm} 
Let $\X$ be a genus $g\geq 2$ superelliptic curve, defined over $\C$, with an extra automorphism, $\s_1, \dots , \s_g$  its dihedral invariants, $F$ the field of moduli, and $K$ its minimal field of definition.  Then, the following are true \\
 
i) The minimal field of definition $K$  is    $K = F ( \sqrt{\Delta_\u})$   

ii)  The equation of $\X$ over $K$ is 
\begin{equation}\label{normal_eq-2}
y^n = A\, x^{\d (s+1)} + A \,  x^{\d s} + \sum_{i=1}^{s-1}  2^{s-i} \, \u_1 \cdot \frac  { \u_s^i \u_i - A \u_{s+1-i}   }    { 2^s \u_1^2 - \u_s^{s+1}  } \cdot x^{\d \cdot i} + 1  
 \end{equation}
where \[ 2^{s+1} A^2 - 2^{s+1} \u_1 A + \u_s^{s+1} = 0. \]  
\end{thm}

\proof
Part i) is an immediate consequences of the above.  To prove part ii) we have to express the coefficients $a_i a_s^i$ of $x^{\d \cdot i}$, $i=2, \dots , s-1$,  in terms of $\u_1, \dots , \u_s$. From the definitions of $\u_i$ we get the following equations:
\[
\left\{
\begin{split}
& \u_1  = a_1^{s+1} + a_s^{s+1}  \\
& \u_s = 2 a_1 a_s \\ 
& \u_i  = a_1^{s+1-i} a_i + a_s^{s+1-i}a_{s+1-i} \\
& \u_{s+1-i} = a_1^{i} a_{s+1-i} + a_s^{i} a_{i}\\
& A=a_s^{s+1} \\
\end{split}
\right.
\]
We multiply both sides in the definition of $\u_i$ by $a_s^i a_1^i = \left( \frac {\u_s} 2 \right)^i$ and have 
\begin{equation}\label{eq_u-i}
 \left( \frac {\u_s} 2 \right)^i \u_i = a_1^{s+1} \cdot a_i \cdot a_s^i + a_s^{s+1} \cdot a_1^i \cdot a_{s+1-i} \end{equation}
From the definition of $\u_{s+1-i}$ we have \[a_1^i a_{s+1-i}= \u_{s+1-i} - a_i a_s^i,\] which we substitute in the Eq.~\eqref{eq_u-i}. 
Hence, 
\[ \boxed{ a_i a_s^i = \frac 1 {a_1^{s+1} - a_s^{s+1} } \, \left( \frac {\u_s^i} {2^i} \u_i - A \, \u_{s+1-i} \right) } \]
Denote by $B:={a_1^{s+1} - a_s^{s+1}}$.  Notice that 

\[ 
\begin{split}
 \left( a_1^{s+1} - a_s^{s+1} \right)  \left( a_1^{s+1} + a_s^{s+1} \right) & =   a_1^{2(s+1)} - a_s^{2(s+1)}  \\
    & = a_1^{2(s+1)} + 2 \left(  a_1 a_s \right)^{s+1} + a_s^{2(s+1)} - 2 \left( a_1 a_s \right)^{s+1} \\
    & = \left(  a_1^{s+1} + a_s^{s+1}\right)^2 - 2 \, \left(  \frac {\u_s} 2 \right)^{s+1} \\
    & =  \u_1^2 - \frac 1 {2^s}     \, \u_s^{s+1}  \\
\end{split}
\]    
Hence,  $ B \u_1 = \u_1^2 - \frac 1 {2^s}     \, \u_s^{s+1}$ 
and
\[ B = \u_1 - \frac 1 {2^s} \, \frac {\u_s^{s+1}} {\u_1}, \]
provided that $\u_1 \neq 0$. 

Hence,
\[ a_i a_s^i = 2^{s-i} \, \u_1 \cdot \frac  { \u_s^i \u_i - A \u_{s+1-i}   }    { 2^s \u_1^2 - \u_s^{s+1}  } \]
as claimed.  This completes the proof. 
\endproof


The natural question is for what values of $\u_1, \dots, \u_s$ is 
\[ \D_\u =  2^{s+1} \left(  2^{s+1}  \u_1^2 - 4 \u_s^{s+1}  \right)         \]  
a complete square in $K$. In this case the field of moduli would be equal to the field of definition.

\bibliographystyle{amsplain}

\end{document}